\documentclass{IEEEtran4PSCC}
\ifCLASSINFOpdf
\else
\fi
\usepackage[T1]{fontenc}
\usepackage[utf8]{inputenc}
\usepackage[cmex10]{amsmath}
\usepackage{textcomp}
\usepackage{upgreek}
\usepackage{hyperref}
\usepackage{cleveref}
\usepackage{commath}
\usepackage{comment}
\usepackage{float}
\usepackage{graphicx}
\usepackage{listings}
\usepackage{tikz}
\usepackage[makeroom]{cancel}
\usepackage{gensymb}
\usepackage{makecell}
\usepackage{booktabs}
\usepackage{dcolumn}
\usepackage{wrapfig}
\usepackage{mathtools}
\usepackage{siunitx}
\usepackage[caption=false,font=footnotesize,position=top]{subfig}
\usepackage[official]{eurosym}
\usepackage{url}

\crefformat{equation}{(#2#1#3)}
\crefrangeformat{equation}{(#3#1#4) to~(#5#2#6)}
\crefmultiformat{equation}{(#2#1#3)}%
{ and~(#2#1#3)}{, (#2#1#3)}{ and~(#2#1#3)}

\renewcommand{\vec}[1]{\boldsymbol{#1}}
\newcommand{\M}{\mathcal{M}}
\newcommand{\T}{\mathcal{T}}
\newcommand{\LO}{\mathcal{L}}
\newcommand{\I}[2]{\mathcal{I}_{#1}^{#2}}
\newcommand{\N}[1]{\mathcal{N}_{#1}}
\newcommand{\K}{\mathcal{K}}
\newcommand{\J}{\mathcal{J}}

\newcommand{\B}{\vec{B}_3(t)}

\usepackage{nomencl}
\makenomenclature
\usepackage{etoolbox}
\renewcommand\nomgroup[1]{%
  \item[\bfseries
  \ifstrequal{#1}{A}{Sets and indices}{%
  \ifstrequal{#1}{B}{Variables}{%
  \ifstrequal{#1}{C}{Parameters}{}}}%
]}

\newcommand{\nomunit}[1]{%
\renewcommand{\nomentryend}{\hspace*{\fill}#1}}

\begin{document}
%
\title{Hydrothermal Scheduling in the Continuous-Time Framework}


 \author{
 \IEEEauthorblockN{
 Christian Øyn Naversen\IEEEauthorrefmark{1}\IEEEauthorrefmark{2},
 Arild Helseth\IEEEauthorrefmark{3},
 Bosong Li\IEEEauthorrefmark{4},
 Masood Parvania\IEEEauthorrefmark{4},
 Hossein Farahmand\IEEEauthorrefmark{2},
 João P. S. Catalão\IEEEauthorrefmark{8}
 }
 
\IEEEauthorblockA{
\IEEEauthorrefmark{2} Department of Electric Power Engineering, Norwegian University of Science and Technology, Trondheim, Norway}
 
\IEEEauthorblockA{
\IEEEauthorrefmark{3} Department of Energy Systems, SINTEF Energy Research, Trondheim, Norway}

\IEEEauthorblockA{
\IEEEauthorrefmark{4} Department of Electrical and Computer Engineering, University of Utah, Salt Lake City, Utah, USA} 

\IEEEauthorblockA{
\IEEEauthorrefmark{8} Faculty of Engineering of the University of Porto and INESC TEC, Porto, Portugal}

}

 
 
 

\maketitle

\begin{abstract}
Continuous-time optimization models have successfully been used to capture the impact of ramping limitations in power systems. In this paper, the continuous-time framework is adapted to model flexible hydropower resources interacting with slow-ramping thermal generators to minimize the hydrothermal system cost of operation. To accurately represent the non-linear hydropower production function with forbidden production zones, binary variables must be used when linearizing the discharge variables and the continuity constraints on individual hydropower units must be relaxed. To demonstrate the performance of the proposed continuous-time hydrothermal model, a small-scale case study of a hydropower area connected to a thermal area through a controllable high-voltage direct current (HVDC) cable is presented. Results show how the flexibility of the hydropower can reduce the need for ramping by thermal units triggered by intermittent renewable power generation. A reduction of 34\% of the structural imbalances in the system is achieved by using the continuous-time model.
\end{abstract}


\thanksto{\noindent \IEEEauthorrefmark{1}Corresponding author: christian.naversen@ntnu.no\\
This work was funded by The Research Council of Norway, Project No. 268014/E20.}

\section{Introduction}
The Norwegian power system is in an interesting state of transition towards tighter integration to the rest of Europe. New high-voltage direct current (HVDC) cable interconnections to Germany and Great Britain are under construction, which increases the potential of cross-zonal trading of both energy and balancing services. Hydropower dominates the Norwegian generation mix and is well suited to provide system balancing services due to its flexibility. A larger share of intermittent renewable generation means that hydropower will play an increasingly important role in providing flexibility to the interconnected North European system in the future. However, propagating the flexibility across HVDC cables is challenging with current practices related to the hourly day-ahead market structure. According to the Norwegian transmission system operator Statnett, changing the HVDC cable flow between areas on an hourly basis has the potential of increasing the structural (or deterministic) imbalances caused by the mismatch in the scheduled hourly production and real-time load \cite{statnett}. In this paper, a modified version of the continuous-time optimization framework is proposed to impose a smooth and continuous flow of power between a hydropower area and a thermal area connected by an HVDC cable.

Continuous-time optimization was originally used to accurately describe the cost of ramping scarcity in thermal systems with large amounts of renewable power generation, such as the power system in California \cite{Parvania2016}. Ramping restrictions can be directly applied to the derivatives of the decision variables when they are allowed to be continuous and smooth functions of time instead of the usual piece-wise constant formulation. The continuous-time formulation relies on limiting the decision variables to be polynomials of degree $r$, which allows the variables to be expressed by the Bernstein polynomials of the same degree. The optimization problem can then be defined in terms of the coefficients of the Bernstein polynomials, which is a mixed-integer linear program (MILP) in the case of the unit commitment problem. The continuous-time framework has lately been expanded in several directions. The existence of a continuous-time marginal price for the economic dispatch problem was proven and calculated in \cite{Parvania2017} for a thermal system. This work was later extended to include energy storage devices in \cite{Khatami2018d}, which has applications in optimal control of charging electric vehicles according to queue theory \cite{Khatami2018b, Khatami2018a} and the scheduling of batteries in balancing markets \cite{Khatami2019}. A stochastic continuous-time model was formulated for unit commitment and reserve scheduling problem in \cite{Hreinsson2018}, with the inclusion of energy storage in \cite{Hreinsson2019} and a method for load estimation and scenario generation in \cite{Khatami2018c}. Applications to other areas within the power system operations field are also emerging, such as the active distribution network model in \cite{Deng2019}.

Hydrothermal scheduling has been an active field of research for decades, which in turn has contributed to the advanced mathematical models used for system and operational planning in hydropower-dominated systems. Good examples of this are the models used in Norway \cite{Fosso1999,Wolfgang2009,Helseth2018} and Brazil \cite{Diniz2018,Maceiral2018}. Previous hydrothermal scheduling models have been based on the standard discrete-time formulation, which assumes piece-wise constant values for time-dependent variables and model input parameters. This paper concerns the novel integration of hydrothermal scheduling and the continuous-time framework. In particular, the integrated modeling of continuous-time operation of complex hydropower cascades poses several new challenges to both hydropower scheduling and continuity constraints. The novel contributions of this paper are outlined as follows:
\begin{itemize}
    \item A continuous-time model including hydropower, thermal generation, and HVDC cables is formulated and studied. To the best of the authors' knowledge, this has not been published previously.
    \item A method for modelling the forbidden production zone of the hydropower production curve in the continuous-time framework is presented. This involves enforcing the continuity constraints on the sum of generated hydropower instead of on the individual plants.
    \item The issue of correct uploading of piece-wise linearized variables in the continuous-time framework is highlighted in relation to the hydropower production function, and a binary variable solution is presented.
\end{itemize}

\Cref{model_section} presents the novel continuous-time model in detail, which is then solved for a two-area system and compared to a discrete-time (hourly) model in \Cref{case_study_section}. Concluding remarks are given in \Cref{conclusion_section}.

\nomenclature[A]{$\mathcal{A}$}{Areas in the system, index $a$}
\nomenclature[A]{$\M$}{Hydropower plants and reservoirs, index $m$}
\nomenclature[A]{$\N{m}$}{Discharge segments in plant $m$, index $n$}
\nomenclature[A]{$\mathcal{I}^{d/b/o}_{m}$}{Reservoirs that discharge/bypass/spill into $m$, index $i$}
\nomenclature[A]{$\T$}{Time intervals, index $h$}
\nomenclature[A]{$\J$}{Thermal generators, index $j$}
\nomenclature[A]{$\K$}{Water value cuts, index $k$}
\nomenclature[A]{$\LO$}{HVDC cables, index $l$}

\nomenclature[B]{$\alpha$}{Future expected system cost \nomunit{[mu]}}
\nomenclature[B]{$q^{d}_{m}(t)$, $\vec{q}^{d}_{mh}$}{Flow through turbine \nomunit{[$\text{m}^3$/s]}}
\nomenclature[B]{$q^{s}_{mn}(t)$, $\vec{q}^{s}_{mnh}$}{Flow through discharge segment \nomunit{[$\text{m}^3$/s]}}
\nomenclature[B]{$q^{b}_{m}(t)$, $\vec{q}^{b}_{mh}$}{Flow through bypass gate \nomunit{[$\text{m}^3$/s]}}
\nomenclature[B]{$q^{o}_{m}(t)$, $\vec{q}^{o}_{mh}$}{Flow through spill gate \nomunit{[$\text{m}^3$/s]}}
\nomenclature[B]{$q^{net}_{m}(t)$, $\vec{q}^{net}_{mh}$}{Net flow into reservoir \nomunit{[$\text{m}^3$/s]}}
\nomenclature[B]{$q^{in}_{m}(t)$, $\vec{q}^{in}_{mh}$}{Total controlled flow into reservoir \nomunit{[$\text{m}^3$/s]}}
\nomenclature[B]{$q^{out}_{m}(t)$, $\vec{q}^{out}_{mh}$}{Total controlled flow out of reservoir \nomunit{[$\text{m}^3$/s]}}
\nomenclature[B]{$q^{rel}_{m}(t)$, $\vec{q}^{rel}_{mh}$}{Total flow released out of reservoir \nomunit{[$\text{m}^3$/s]}}
\nomenclature[B]{$v_{m}(t)$}{Reservoir volume \nomunit{[$\text{m}^3$]}}
\nomenclature[B]{$p_{m}(t)$, $\vec{p}_{mh}$}{Generated hydropower \nomunit{[MW]}}
\nomenclature[B]{$g_{j}(t)$, $\vec{g}_{jh}$}{Generated thermal power \nomunit{[MW]}}
\nomenclature[B]{$f_{l}(t)$, $\vec{f}_{lh}$}{Flow on HVDC cable \nomunit{[MW]}}
\nomenclature[B]{$s_{j}^{\uparrow/\downarrow}(t)$}{Startup/shutdown of thermal generator \nomunit{[MW]}}
\nomenclature[B]{$\bar{s}_{m}^{\uparrow/\downarrow}(t)$}{Startup/shutdown of hydropower plant \nomunit{[MW]}}
\nomenclature[B]{$w_{mn}(t)$}{Discharge segment commitment decision}
\nomenclature[B]{$u_{j}(t)$}{State of thermal unit, on/off}
\nomenclature[B]{$z_{m}(t)$}{State of hydropower unit, on/off}

\nomenclature[C]{$WV_{mk}$}{Water value cut coefficient \nomunit{[mu/$\text{m}^3$]}}
\nomenclature[C]{$D_{k}$}{Water value cut constant \nomunit{[mu]}}
\nomenclature[C]{$C^{b}$}{Penalty for bypassing water \nomunit{[mu/$\text{m}^3$]}}
\nomenclature[C]{$C^{o}$}{Penalty for spilling water \nomunit{[mu/$\text{m}^3$]}}
\nomenclature[C]{$C_j$}{Marginal cost of thermal generator \nomunit{[mu/MW]}}
\nomenclature[C]{$C_j^{\uparrow/\downarrow}$}{Thermal unit startup/shutdown cost  \nomunit{[mu]}}
\nomenclature[C]{$Q^{s}_{mn}$}{Maximal flow through turbine segment \nomunit{[$\text{m}^3$/s]}}
\nomenclature[C]{$Q^{d}_{m}$}{Maximal flow through turbine \nomunit{[$\text{m}^3$/s]}}
\nomenclature[C]{$Q^{b}_{m}$}{Maximal flow through bypass gate \nomunit{[$\text{m}^3$/s]}}
\nomenclature[C]{$I_{m}(t)$}{Natural inflow into reservoir \nomunit{[$\text{m}^3$/s]}}
\nomenclature[C]{$I^{u}_{m}(t)$}{Natural inflow from creek intakes \nomunit{[$\text{m}^3$/s]}}
\nomenclature[C]{$V^{0}_{m}$}{Initial reservoir volume \nomunit{[$\text{m}^3$]}}
\nomenclature[C]{$V_{m}$}{Maximal reservoir capacity \nomunit{[$\text{m}^3$]}}
\nomenclature[C]{$\eta_{mn}$}{Energy conversion factor \nomunit{[MWs/$\text{m}^3$]}}
\nomenclature[C]{$P_{m}^{max/min}$}{Maximal/minimal hydropower capacity \nomunit{[MW]}}
\nomenclature[C]{$G_{j}^{max/min}$}{Maximal/minimal thermal capacity \nomunit{[MW]}}
\nomenclature[C]{$L_a(t)$}{Net area load \nomunit{[MW]}}
\nomenclature[C]{$\delta_{h}$}{Length of time interval \nomunit{[s]}}
\nomenclature[C]{$N$}{Number of time intervals in $\T$}
\nomenclature[C]{$G_{la}$}{Line flow direction coefficient}
\nomenclature[C]{$F_l^{max}$}{Maximal flow limit on HVDC cable \nomunit{[MW]}}
\nomenclature[C]{$R_l^{u/d}$}{Ramping limits of HVDC cable flow \nomunit{[MW/s]}}
\nomenclature[C]{$R_j^{u/d}$}{Ramping limits of running thermal unit \nomunit{[MW/s]}}
\nomenclature[C]{$R_j^{\uparrow/\downarrow}$}{Thermal ramping gain for starts/stops \nomunit{[MW/s]}}

\printnomenclature[1.95cm]

\section{Model}\label{model_section}
\subsection{Fundamentals of a continuous-time model}
The core idea of the continuous-time framework is to represent time-dependent input and decision variables as polynomials of time instead of piece-wise constant functions. This increases the complexity of the model formulation, but sub-hourly effects and constraints related to derivatives with respect to time are easily captured. The motivation behind the original continuous-time unit commitment model in \cite{Parvania2016} was precisely to incorporate the impact of ramping scarcity into the market clearing. The time-dependent decision variables in the typical continuous-time optimization framework are defined through the Bernstein polynomials of degree $r$, $\vec{B}_r(t)$, which form a basis for any polynomials of at most degree $r$ on the time interval $[0,1]$. By splitting the time horizon of the model into $N$ intervals $h\in\T$ of length $\delta_h$, the time-dependent decision variables can be expressed as polynomials of the form

\begin{equation}\label{bernstein_exp}
x(t) = \sum_{h\in\T} \vec{x}^T_h \cdot \vec{B}_r(\tau_h)\Pi(\tau_h),
\end{equation}
where $\tau_h$ and $\Pi(\tau_h)$ are defined as follows:
\begin{equation}\label{tau_trans}
\tau_h = \frac{1}{\delta_h}\Bigg(t-\smashoperator{\sum_{i<h}}\delta_i\Bigg) \hspace{60pt} \forall h\in\T,
\end{equation}
\begin{equation}
\Pi(\tau_h) = \left\{
        \begin{array}{ll}
            1, & \quad 0\leq \tau_h \leq 1 \\
            0, & \quad \text{otherwise}
        \end{array}
    \right.\qquad \forall h\in\T.
\end{equation}

The vectors $\vec{x}_h$ contain the $r+1$ coefficients of the Bernstein polynomials in each time interval, which become the decision variables of the continuous-time model. It is necessary to use the scaled time $\tau_h$ and the operator $\Pi$ to project the Bernstein polynomials into the correct time interval while maintaining their property as basis functions. One of the main reason for using Bernstein polynomials is the convex hull property, which makes it possible to impose inequality constraints on $x(t)$ for all times $t$ by directly bounding the coefficients $\vec{x}_h$ \cite{Parvania2016}. This paper uses Bernstein polynomials of degree 3 as the basis:

\begin{equation}\label{bernstein_def}
\B = \begin{bmatrix}(1-t)^3,& 3t(1-t)^2,& 3t^2(1-t),& t^3 \end{bmatrix}^T.
\end{equation} 

This is a popular choice in the literature, as it keeps the size of the model reasonable without sacrificing the ability to model complex time dependencies. Another advantage is the linear relationship to the cubic Hermite splines $\vec{H}(t)$, which can be used as an equivalent basis:


\begin{equation} \label{hermite_relation}
\vec{H}(t) = 
\begin{bmatrix}
1 & 1 & 0 & 0 \\
0 & \frac{1}{3} & 0 & 0 \\
0 & 0 & 1 & 1 \\
0 & 0 & -\frac{1}{3} & 0 \\
\end{bmatrix} \cdot \B \equiv \vec{\mathbf{W}}\cdot\B.
\end{equation}

The coefficients of the Hermite splines have a physical interpretation as the value of $x(t)$ and its derivative $\dot{x}(t)$ at the start and end of the time interval $h$:

\begin{equation}\label{continuity_prop}
\vec{x}^{\vec{H}}_h = (\vec{\mathbf{W}}^{-1})^T\cdot\vec{x}_h = \Big[x_h^{start},\hspace{2pt} \dot{x}_h^{start},\hspace{2pt} x_h^{end},\hspace{2pt} \dot{x}_h^{end}\Big]^T.
\end{equation}

This interpretation is useful for expressing the continuity of $x(t)$ across the time intervals $h$. The reader is referred to \cite{Parvania2016} for a more detailed introduction to the continuous-time formulation with further references to the properties of the Bernstein polynomials mentioned in this section.

\subsection{Objective function}
The objective of the proposed hydrothermal model is to minimize the future expected cost of the system, the penalties for bypassing and spilling water, and the operational, startup and shutdown costs of the thermal generators:

\begin{align}
Z &= \alpha + \sum_{m\in\M}\smashoperator{\int_{0}^{t^{end}}}\left(C^{b}q_{m}^{b}(t) + C^{o}q_{m}^{o}(t)\right)dt\nonumber\\
&+\sum_{j\in\J}\smashoperator{\int_0^{t^{end}}}C_j g_j(t)dt + \sum_{j\in\J}\sum_{h\in\T}\left(C_j^{\uparrow}s^{\uparrow}_{jh} +  C_j^{\downarrow}s^{\downarrow}_{jh}\right).\label{obj_integral}
\end{align}

Note that startup and shutdown cost are assumed to be negligible for the hydropower plants. The definite integral of the Bernstein polynomials of the third degree is $\int_0^1 \vec{B}_3(t)dt = \frac{1}{4}\vec{1}$, which simplifies the integrals in \cref{obj_integral} to the sums

\begin{align}
Z &= \alpha + \frac{1}{4}\sum_{m\in\M}\sum_{h\in\T}\delta_h\vec{1}^T\cdot\left(C^{b}\vec{q}_{mh}^{b} + C^{o}\vec{q}_{mh}^{o}\right)\nonumber\\
&+\sum_{j\in\J}\sum_{h\in\T} \left(\frac{1}{4}\delta_hC_j\vec{1}^T\cdot \vec{g}_{jh} + C_j^{\uparrow}s^{\uparrow}_{jh} +  C_j^{\downarrow}s^{\downarrow}_{jh}\right).\label{obj}
\end{align}

As this paper focuses on modelling hydropower generation in the continuous-time framework, a simplified linear formulation of the thermal generation cost function is used in \cref{obj}. More advanced modeling of quadratic and piece-wise linear cost functions in continuous-time unit commitment are available in the literature \cite{Parvania2016, Hreinsson2019}, and their integration in the model proposed in this paper is straightforward. 

\subsection{Hydropower topology constraints}
The cascaded topology constraints dictate how water moves between the reservoirs. These constraints are equality constraints, see for instance \cite{Helseth2018}, which means that equating the polynomial coefficients are sufficient to satisfy them in the continuous-time framework. The convex hull property of the Bernstein polynomials and the fact that $\vec{1}^T\cdot\vec{B}_3(t) = 1$ is used to enforce the physical bounds on the variables:

\begin{align}
&\vec{q}_{mh}^{net} = \vec{I}_{mh} + \vec{q}_{mh}^{in} - \vec{q}_{mh}^{out} \hspace{52pt} \forall m,h\in\M,\T \label{q_net}\\
&\vec{q}_{mh}^{out} = \vec{q}_{mh}^{rel} +\vec{q}_{mh}^{o} \hspace{82pt} \forall m,h\in\M,\T  \label{q_out}\\
&\vec{q}_{mh}^{in} =  \sum_{i\in\I{m}{d}}\vec{q}_{ih}^{d} + \sum_{i\in\I{m}{b}}\vec{q}_{ih}^{b} + \sum_{i\in\I{m}{o}}\vec{q}_{ih}^{o} \hspace{1pt} \forall m,h\in\M,\T  \label{q_in}\\
&\vec{q}_{mh}^{rel} = \vec{q}_{mh}^{d} + \vec{q}_{mh}^{b} - \vec{I}_{mh}^{u} \hspace{50pt} \forall m,h\in\M,\T  \label{q_rel}\\
&\vec{0}\leq \vec{q}_{mh}^{d} \leq Q_{m}^{d}\vec{1} \hspace{88pt} \forall m,h \in \M,\T \label{dis_bounds}\\
&\vec{0}\leq \vec{q}_{mh}^{b} \leq Q_{m}^{b}\vec{1} \hspace{88pt} \forall m,h \in \M,\T \label{bypass_bounds}\\
&\vec{0}\leq \vec{q}_{mh}^{o}  \hspace{124pt} \forall m,h \in \M,\T \label{spill_bounds}\\
&\vec{0}\leq \vec{q}_{mh}^{rel} \hspace{124pt} \forall m,h \in \M,\T. \label{rel_bounds}
\end{align}

\begin{figure}[ht]
    \centering 
    \includegraphics[width=0.5\columnwidth]{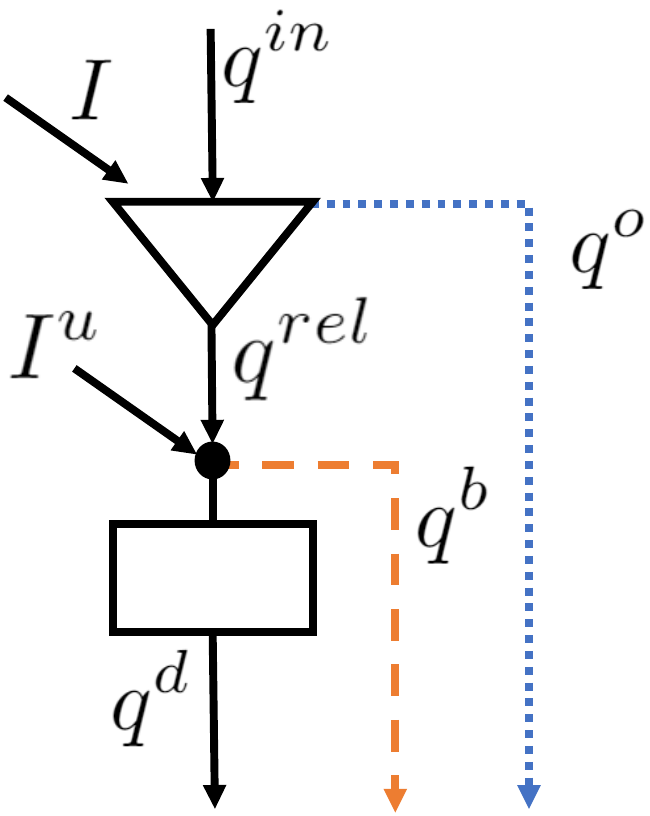}
\caption{Depiction of the different waterways for discharging, bypassing and spilling water between reservoirs. All waterways may lead to different downstream reservoirs or out of the system. Natural inflow enters the system in two different ways, either into the reservoir (triangle shape) or directly into the main tunnel of the plant (rectangle shape).}
\label{fig:module}
\end{figure}

There are three waterways that connect reservoirs; discharge through the turbine, the bypass gate and the spill gate. \Cref{fig:module} shows the relationship between the different waterways in addition to where natural inflow enters the system.

\subsection{Volume constraints}

The rate of change in the reservoir content is described by the differential equation:

\begin{equation}
\frac{dv_{m}(t)}{dt} = q_{m}^{net}(t) \qquad \forall m\in\M.
\end{equation}

The integral of Bernstein polynomials of degree 3 can be expressed using Bernstein polynomials of degree 4 using a linear mapping matrix \cite{Khatami2018d,Hreinsson2019}:

\begin{equation}
\smashoperator{\mathop{\int\B dt}} = \frac{1}{4}
\begin{bmatrix}
0 & 1 & 1 & 1 & 1\\
0 & 0 & 1 & 1 & 1\\
0 & 0 & 0 & 1 & 1\\
0 & 0 & 0 & 0 & 1\\
\end{bmatrix} \vec{B}_4(t) \equiv \vec{\mathbf{N}}\cdot\vec{B}_4(t)
\end{equation}
which is further utilized to show the volume increase within a time interval $h$ as follows:

\begin{align}
&v_{m}(t) -v_{mh} = \smashoperator{\int_{t(h)}^{t}} \frac{dv(t')}{dt'}dt' = (\vec{q}_{mh}^{net})^T\cdot \smashoperator{\int_{t(h)}^{t}}\vec{B}_3(\tau_h')dt\nonumber\\
&= \delta_h(\vec{q}_{mh}^{net})^T\cdot\vec{\mathbf{N}}\cdot \vec{B}_4(\tau_h).\label{vol_calc}
\end{align}

Here, $t(h) = \sum_{i<h}\delta_i$ and the fact that $\vec{\mathbf{N}}\cdot\vec{B}_4(0) = \vec{0}$ was used. Note that the volume variables $v_{mh}$ denotes the volume at the start of interval $h$. Based on equation \cref{vol_calc}, the following volume balance constraints can be added to the optimization problem:

\begin{align}
&v_{m0} = V_{m}^0 \hspace{117pt} \forall m\in\M \label{vol_init}\\
&v_{m,h+1}-v_{mh} = \frac{1}{4}\delta_h\vec{1}^T{\cdot}\vec{q}_{mh}^{net} \hspace{38pt} \forall m,h\in\M,\T \label{vol_bal}\\
&\vec{0} \leq v_{mh}\vec{1} + \delta_h\vec{\mathbf{N}}^T\cdot\vec{q}_{mh}^{net} \leq V_{m}\vec{1} \hspace{23pt}  \forall m,h\in\M,\T. \label{vol_bounds}
\end{align}

Constraint \cref{vol_init} sets the initial volume of each reservoir and \cref{vol_bal} calculates the volume change from one time interval to the next by inserting $\vec{\mathbf{N}}\cdot\vec{B}_4(1) = \frac{1}{4}\vec{1}$. Constraint \cref{vol_bounds} uses the convex hull property to bound the volume within the limits of the reservoir for all times $t$.

\subsection{Future cost bounds}
The future expected cost of the system is represented by a set of Benders cuts created by a hydrothermal long-term model such as \cite{Helseth2018}. The expected future cost depends on the state of all hydropower reservoirs in the system at the end of the last time interval $N$:

\begin{equation}
\alpha \geq \sum_{m\in\M}WV_{mk}v_{m,N+1} + D_{k} \qquad \forall k \in \K. \label{cuts}    
\end{equation}

\subsection{Hydropower production} \label{prod_function_section}

The conversion from discharge through the turbine to generated power is a non-linear function which depends on the effective plant head and the efficiency curves of the turbine and generator \cite{Kong2020}. By assuming a constant head for the planning horizon, the hydropower production function can be approximated as a single piece-wise linear curve, where the discharge variable is split into $n\in\N{m}$ segments with constant gradient $\eta_n$. In an discrete-time model, the discharge segments will usually be uploaded in the correct order as long as the gradient is decreasing for increasing segment number. The exception is extreme situations where it is beneficial to dump as much water as possible while limiting the power produced, which can be the case in high inflow and low load scenarios. A similar effect of incorrect uploading of discharge segments has been observed in this work when the continuous-time framework was implemented. Segments with high efficiency are still favoured but there is no guarantee that segment $n$ is at its maximal capacity for all times that segment $n+1$ is being used. The model will often start using the next segment too early to be able to fulfill the continuous-time power balance described in \Cref{power_bal_section}. To remedy this problem, binary variables $w_{mn}(t) = \sum_{h\in\T} w_{mnh} \vec{1}^T\cdot\vec{B}_3(\tau_h)\Pi(\tau_h)$ are used in this work to force the segments to be fully utilized before the next segment can be used:

\begin{align}
&\vec{q}_{mh}^{d} = \sum_{n\in\N{m}}\vec{q}_{mnh}^{s} \hspace{50pt} \forall m,h \in \M,\T \label{dis_sum}\\
&\vec{p}_{mh} = \sum_{n\in\N{m}}\eta_{mn}\vec{q}_{mnh}^{s} \hspace{32pt} \forall m,h\in \M,\T \label{pq}\\
& Q_{mn}^{s}w_{mnh}\vec{1} \leq\vec{q}_{mnh}^{s} \leq Q_{mn}^{s}\vec{1} \hspace{5pt} \forall m,h,n \in \M,\T,\N{m} \label{dis_seg_bound_1}\\
&\vec{q}_{mnh}^{s} \leq Q_{mn}^{s}w_{m,n-1,h}\vec{1} \hspace{8pt}  \forall m,h,n \in \M,\T,\N{m}\backslash\{0\}. \label{dis_seg_bound_2}
\end{align}

This modelling choice of the hydropower production function has the unfortunate effect of introducing additional binary variables into the model but also enables the use of non-concave linearizations of the hydropower production function. It is also possible to incorporate forbidden production regions within the operating range of the turbine by modifying \cref{dis_seg_bound_1} to $\vec{q}_{mnh}^s = Q_{mn}^{s}w_{mnh}\vec{1}$ for the segment representing the forbidden region. 

\subsection{Power balance and HVDC power flow}\label{power_bal_section}

The power balance constraints must be satisfied in each node of the system. In this work, each node represents a larger market area assuming no internal power flow limits. The areas are connected with HVDC cables where the flow can be controlled by the system operator. The power balance constraints are formulated as

\begin{equation}\label{power_bal}
\smashoperator{\sum_{m\in\M_{a}}}\vec{p}_{mh} + \smashoperator{\sum_{j\in\J_{a}}}\vec{g}_{jh} - \smashoperator{\sum_{l\in\LO}} G_{la}\vec{f}_{lh} = \vec{L}_{ah} \hspace{10pt} \forall a,h \in \mathcal{A},\T.
\end{equation}

The coefficient $G_{la}$ dictates the positive and negative direction of flow on each cable $l\in\LO$ by taking the values $\pm1$, or zero if cable $l$ is not connected to area $a$. $\M_{a}$ and $\J_{a}$ are the sets of hydropower and thermal units located in area $a$, respectively. The flow on the HVDC cables is constrained by maximal flow limits

\begin{equation}
-F_l^{max}\vec{1} \leq \vec{f}_{lh} \leq F_l^{max}\vec{1}    \qquad \forall l,h \in \LO,\T,
\end{equation}
and additional limitations on the change of flow is imposed on the derivative $\dot{f}_l(t)$ to stay within the specified HVDC cable ramping limits used in the Nordic system \cite{nordpool}. By using the following property of the Bernstein polynomials,

\begin{equation}
\frac{d\B}{dt} = 3
\begin{bmatrix}
-1 & 0 & 0\\
1 & -1 & 0\\
0 & 1 & -1\\
0 & 0 & 1 
\end{bmatrix} \cdot \vec{B}_2(t) \equiv \vec{\mathbf{K}}\cdot\vec{B}_2(t),
\end{equation}
the minimum and maximum ramping limits can be expressed as:

\begin{equation}\label{hvdc_ramping}
-R_{l}^{d}\vec{1}^T \leq \frac{1}{\delta_h} \vec{f}_{lh}^T\cdot\vec{\mathbf{K}} \leq R_{l}^{u}\vec{1}^T \qquad \forall l,h \in \LO,\T.
\end{equation}    

\subsection{Thermal generation constraints}

The thermal generators are subject to unit commitment decisions which signify if a generator is offline or producing between the minimal and maximal production limits. The thermal unit commitment constraints are modelled by the use of the binary decision variables $u_j(t)$:
\begin{align}
&G_{j}^{min} \vec{u}_{jh} \leq \vec{g}_{jh} \leq G_{j}^{max} \vec{u}_{jh} \hspace{26pt} \forall j,h\in\J,\T \label{thermal_uc}\\
&\vec{u}_{jh} = \begin{bmatrix}u_{jh},u_{jh},u_{j,h+1},u_{j,h+1}\end{bmatrix}^T  \hspace{5pt}\forall j,h\in\J,\T\backslash\{N\}\label{thermal_binary_vec}\\
&\vec{u}_{jN} = u_{jN}\vec{1}  \hspace{92pt}\forall j\in\J \label{thermal_binary_last_vec}\\
&s_{jh}^{\uparrow} - s_{jh}^{\downarrow} = u_{j,h+1} - u_{jh} \hspace{40pt} \forall j,h \in \J,\T \backslash\{N\}\label{thermal_startup}\\
&s_{jh}^{\uparrow} + s_{jh}^{\downarrow} \leq 1 \hspace{79pt}  \quad \forall j,h\in\J,\T \label{thermal_unique_startstop}\\
&u_{jh}, s_{jh}^{\uparrow/\downarrow} \in \{0,1\} \hspace{74pt} \forall j,h\in\J,\T. \label{thermal_binary_def}
\end{align}

The constraints closely follow the implementation used in \cite{Parvania2016} and \cite{Hreinsson2018}, which are in turn adapted from the standard discrete-time unit commitment formulation found in for instance \cite{Carrion2006}. The choice of the commitment decision vector in \cref{thermal_binary_vec,thermal_binary_last_vec} allows the thermal generator to use time interval $h$ to ramp up from zero to above $G^{min}$, or conversely ramp down production to zero. The smooth transition is necessary for the continuity constraints that will be applied to the thermal production variables in \Cref{continuity_section}. Constraint \cref{thermal_startup} captures the startups and shutdowns of the generators, which are accounted for in the objective function \cref{obj}. The up and down ramping constraints of thermal generators, taking into account the startup and shutdown ramp limitations, are modeled as follows:

\begin{align}
&\frac{1}{\delta_h} \vec{g}_{jh}^T\cdot\vec{\mathbf{K}} \leq \left(R_{j}^{u} + R_{j}^{\uparrow}s_{jh}^{\uparrow}\right)\vec{1}^T \hspace{18pt} \forall j,h \in \J,\T\label{thermal_ramping_up}\\
&\frac{1}{\delta_h} \vec{g}_{jh}^T\cdot\vec{\mathbf{K}} \geq -\left(R_{j}^{d} + R_{j}^{\downarrow}s_{jh}^{\downarrow}\right)\vec{1}^T \hspace{10pt}\forall j,h \in \J,\T.\label{thermal_ramping_down}
\end{align}

The minimum up and down time constraints of thermal generation is not considered in this paper, and the readers are referred to our previous works for details on modeling these constraints in the continuous-time unit commitment model \cite{Parvania2016}.

\subsection{Hydropower unit commitment}\label{hydro_uc_section}
Due to operating characteristics such as mechanical vibration or loss of efficiency, hydropower turbines usually have one or several forbidden production regions depending on the turbine type. It is important to model these regions when looking at short-term scheduling of a hydropower system to have an accurate representation of the operating range of the hydropower plants. The unit commitment constraints of the hydropower plants in the continuous-time optimization model must account for the forbidden production region so that the flexibility of the plant is not overestimated. The hydropower unit commitment decisions $z_m(t)$ are used to model this in the following way:
\begin{align}
&P_{m}^{min}\vec{z}_{mh} \leq \vec{p}_{mh} \leq P_{m}^{max}\vec{z}_{mh} \hspace{9pt} \forall m,h \in \M,\T \label{hydro_uc}\\
&\vec{z}_{mh} = z_{mh}\vec{1} \hspace{83pt} \forall m,h \in \M,\T \label{hydro_binary_vec}\\
&\bar{s}_{mh}^{\uparrow} - \bar{s}_{mh}^{\downarrow} = z_{m,h+1} - z_{mh} \hspace{20pt} \forall m,h \in \M,\T\backslash\{N\} \label{hydro_startup}\\
&\bar{s}_{mh}^{\uparrow} + \bar{s}_{mh}^{\downarrow} \leq 1 \hspace{74pt} \forall m,h \in \M,\T \label{hydro_unique_startstop}\\
&z_{mh},\bar{s}_{mh}^{\uparrow/\downarrow} \in \{0,1\} \hspace{64pt} \forall m,h \in \M,\T. \label{hydro_binary_def}
\end{align}

In contrast to the choice of the thermal unit commitment vector in \cref{thermal_binary_vec}, the formulation in \cref{hydro_binary_vec} forces the hydropower unit commitment decisions to be constant for the whole time interval so that the production is never between 0 and $P^{min}$. However, this formulation is in opposition to the normal continuous-time formulation, as discontinuous jumps in power production must be allowed. If not, the hydropower plants will be unable to start and stop at all. These issues are addressed in \Cref{continuity_section}.

\subsection{Continuity constraints}\label{continuity_section}
The standard continuous-time optimization framework builds on the $C^1$ continuity of all decision variables $x(t)$. This requires both the value $x(t)$ and the value of the derivative $\dot{x}(t)$ to be continuous over the change of time intervals $h\in\T$. Such constraints are enforced by using the relationship between the Bernstein polynomials and the cubic spline functions, shown in \cref{hermite_relation}. The interpretation of the coefficients of $\vec{H}(t)$ described in \cref{continuity_prop} simplifies the implementation of the $C^1$ continuity constraints. By labelling the components of the vector $\vec{x}$ as $\vec{x}[i]$ for $i\in\{0,1,2,3\}$, the continuity constraints become:

\begin{align}
&\vec{x}_{h}^{\vec{H}}[2] = \vec{x}_{h+1}^{\vec{H}}[0] \qquad \forall h\in\T\backslash\{N\} \label{value_cont}\\
&\vec{x}_{h}^{\vec{H}}[3] = \vec{x}_{h+1}^{\vec{H}}[1] \qquad \forall h\in\T\backslash\{N\}. \label{deriv_cont}
\end{align}

These constraints are applied to the thermal generation and HVDC flow variables:

\begin{align}
&\vec{g}_{jh}^{\vec{H}}[2] = \vec{g}_{j,h+1}^{\vec{H}}[0] \hspace{5pt}\qquad \forall j,h\in\J,\T\backslash\{N\}\\
&\vec{g}_{jh}^{\vec{H}}[3] = \vec{g}_{j,h+1}^{\vec{H}}[1] \hspace{5pt}\qquad \forall j,h\in\J,\T\backslash\{N\}\\
&\vec{f}_{lh}^{\vec{H}}[2] = \vec{f}_{l,h+1}^{\vec{H}}[0] \hspace{5pt}\qquad \forall l,h\in\LO,\T\backslash\{N\}\\
&\vec{f}_{lh}^{\vec{H}}[3] = \vec{f}_{l,h+1}^{\vec{H}}[1] \hspace{5pt}\qquad \forall l,h\in\LO,\T\backslash\{N\}.
\end{align}

As mentioned in \Cref{hydro_uc_section}, discontinuous jumps in power production are required to model the forbidden production region of hydropower plants. Therefore, enforcing the $C^1$ continuity constraints on the variables related to the hydropower production is not possible. In addition, requiring continuous derivatives for water flow and hydropower production is strict when $\delta_h$ is longer than a few minutes. To avoid conservative solutions underestimating the ramping capabilities of hydropower, \cref{deriv_cont} is not implemented for any variable related to hydropower. The bypass and overflow variables are $C^0$ continuous:

\begin{align}
&\vec{q}_{mh}^{b,\vec{H}}[2] = \vec{q}_{m,h+1}^{b,\vec{H}}[0] \qquad \forall m,h\in\M,\T\backslash\{N\}\\
&\vec{q}_{mh}^{o,\vec{H}}[2] = \vec{q}_{m,h+1}^{o,\vec{H}}[0] \qquad \forall m,h\in\M,\T\backslash\{N\},
\end{align}
and the reservoir volume continuity is already secured by \cref{vol_bal}. The hydropower production is forced to be $C^0$ continuous unless a startup or shutdown happens in the time interval. This is modelled by replacing \cref{value_cont} by the following two inequalities:

\begin{align}
&\vec{p}_{mh}^{\vec{H}}[2] - \vec{p}_{m,h+1}^{\vec{H}}[0] \leq P_{m}^{max}\bar{s}_{mh}^{\downarrow} \hspace{3pt}  \forall m,h\in\M,\T \backslash\{N\} \label{hydro_cont_1}\\
&\vec{p}_{m,h+1}^{\vec{H}}[0] - \vec{p}_{mh}^{\vec{H}}[2] \leq P_{m}^{max}\bar{s}_{mh}^{\uparrow} \hspace{3pt} \forall m,h\in\M,\T\backslash\{N\} \label{hydro_cont_2}
\end{align}
which is consistent with the unit commitment constraints imposed in \crefrange{hydro_uc}{hydro_binary_def}. Note that this relaxation produces a more constrained problem, as production in the forbidden region is impossible. Due to the connection between production and discharge in \cref{pq}, the discharge variables $\vec{q}^s$ must also be allowed to have discontinuous jumps. However, the binary definitions of the discharge bounds in \cref{dis_seg_bound_1,dis_seg_bound_2} take care of continuity when the hydropower plant is producing, so there is no need to apply any further constraints to the discharge variables. The continuity properties of the derived flow variables $\vec{q}^{net}$, $\vec{q}^{out}$, $\vec{q}^{in}$, and $\vec{q}^{rel}$ are also implicitly accounted for through \crefrange{q_net}{q_rel}. 

It is important to note that even though the individual hydropower plants may have discontinuous jumps and discontinuous derivatives in the power production curve between time intervals, their sum is still forced to be $C^1$ continuous through the power balance constraint \cref{power_bal} since all other quantities in the equation are $C^1$ continuous. The $C^1$ continuity constraints of the flexible hydropower have effectively been lifted from the individual plant to the sum on an area level. The hydropower model formulation presented in this paper can be seen as an approximation of a fully $C^1$ continuous model where short time intervals have been inserted around every major time interval shift. By forcing the hydropower plants to only start or stop in these short intervals, an accurate production profile spending minimal time in the forbidden production zone would be achieved. By letting the length of short intervals go to zero, the partially $C^0$ continuous hydropower formulation used in this paper is recovered. Therefore, the alterations made to the continuity constraints for the hydropower-related variables will not drastically impact the operation of the hydropower, as long as $\delta_h$ is long compared to the time it takes to ramp up and down a hydropower plant, which is usually only a few minutes.

\section{Case study}\label{case_study_section}
A small scale case study with two areas connected by a single HVDC cable is presented in this section. The continuous-time model proposed in \Cref{model_section} and an analogous discrete-time hourly model are both solved to show and compare the interaction between fast and slow ramping components in the system. Both models have been implemented in Pyomo and solved with CPLEX 12.8. One area contains only hydropower, while the other only contains thermal generation. The hydropower topology is based on a real Norwegian water course consisting of 12 reservoirs and plants which is described in more detail in \cite{Naversen2020}, and the future expected cost of the hydropower system is calculated based on the long-term model described in \cite{Helseth2018}. The inflow is considered piece-wise constant within each hour in the entire hydropower area, which has a total hydropower production capacity of 537 MW. The thermal area contains four thermal generators with a total of 256 MW of production capacity and varying ramping capabilities and marginal, startup and shutdown costs. The areas are connected by an HVDC cable with a flow limit of 50 MW in either direction. The ramping limitations of the cable are based on the current practice of how fast the flow on an HVDC cable can be changed in the Nordic market, which is 600 MW/h \cite{nordpool}. The flow change is performed in a 20 min window around hourly shifts, which gives an effective ramping rate of 30 MW/min or 1800 MW/h \cite{statnett}. The time horizon is set to 24 hours with hourly time intervals for both the hourly and the continuous-time model.

\begin{figure}[ht]
    \centering 
    \includegraphics[width=\columnwidth]{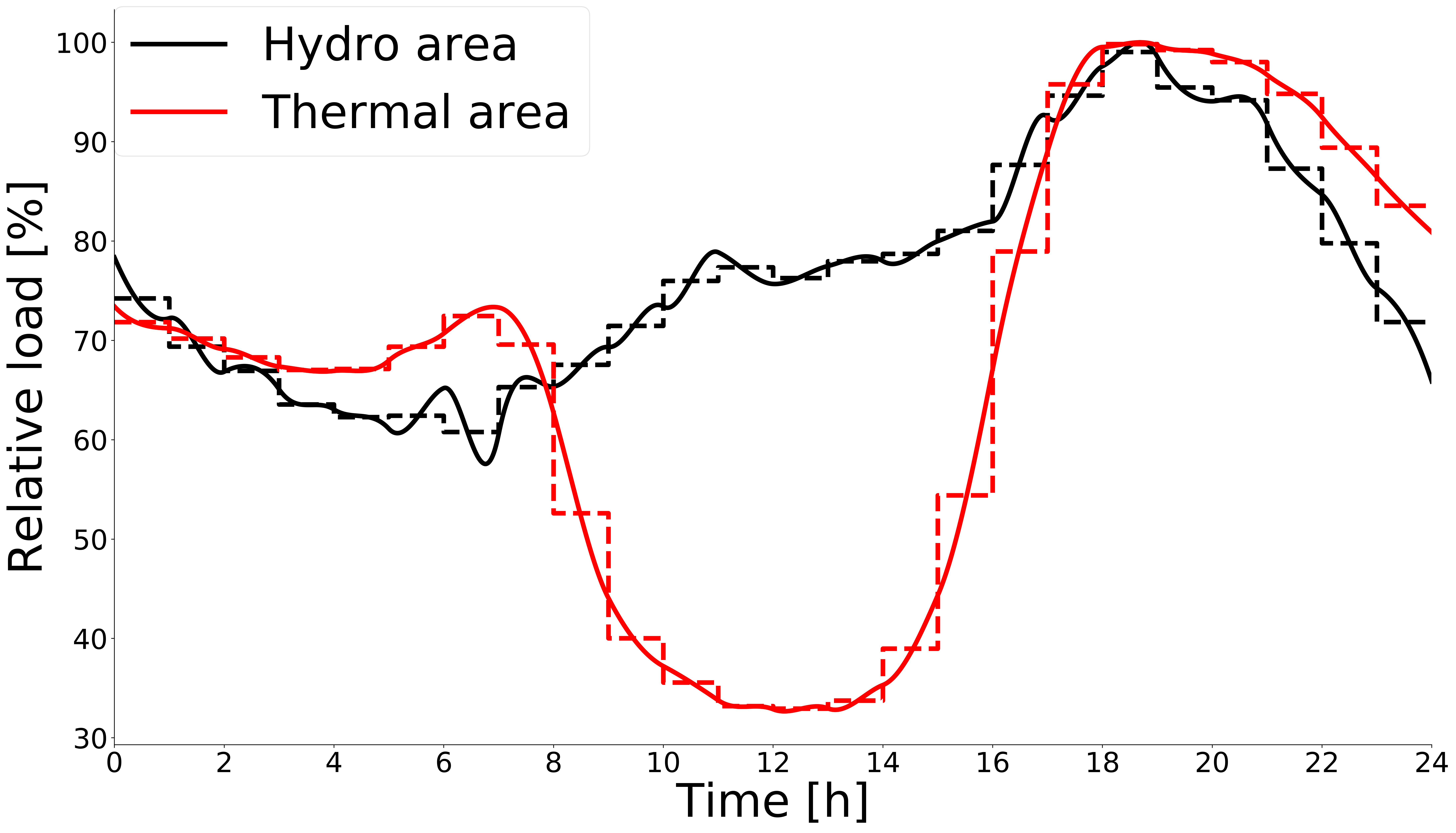}
\caption{The continuous-time load profiles of the thermal and hydropower areas are shown together with the hourly constant load approximations (solid and dashed lines, respectively). The profiles have been scaled by the value of the peak load.}
\label{fig:load}
\end{figure}

The scaled net load profiles for the two areas are shown in \Cref{fig:load}. The peak net load value in each area is used as a scale in the figure, which is 450 MW and 160 MW in the hydropower and thermal areas, respectively. The net load profiles are based on measured data from NYISO and CAISO from 1/1-2019, available at \cite{nyiso,caiso} with a 5-minute resolution. The CAISO net load has been used for the thermal area, which experiences significant ramping events in the morning and afternoon as solar plants start and stop producing power. The continuous-time load was calculated from the raw data by a standard least-squares error fit to the Bernstein polynomials, while the hourly load is the average load for each hour. The structural imbalances in both areas go down in the continuous-time model compared to the hourly model, with a reduction of 34\% on system level. This represents 97 MWh of saved balancing energy, which is 0.9\% of the total daily net system load. The reduction of imbalances is higher in the thermal area (87\%) than the hydropower area (20\%) because $\B$ provides a better fit to the CAISO load data.

The size and solution times of the models are listed in \Cref{table:comparison}, which shows the initial model size and the reduced size after CPLEX performs an automatic presolve routine. The number of continuous and binary variables and constraints are considerably higher in the continuous-time model compared to the hourly model, also after the presolve. The larger model size of the continuous-time model results in a longer solution time on a standard office laptop, i7-7600 CPU at 2.8 GHz with 4 cores, though solution time in MIP models can vary greatly based on the parameter settings given to the solver. A small relative MIP gap of 0.28\% was reached in 60 seconds for the continuous-time model, but solving it to zero gap like the hourly model takes about 10 hours on a server with 36 cores. Upon investigation, it is clear that the hydropower production continuity constraints, \cref{hydro_cont_1,hydro_cont_2}, are the complicating constraints. If these constraints are removed, which means the hydropower production variables are discontinuous over the interval changes, the continuous-time model can be solved to zero MIP gap in 22 seconds. This is a trade-off between realistic physical modelling and tractability that should be considered when solving larger systems.

\begin{table}[hb]
\centering
\caption{Model size comparison of the continuous-time and hourly models. The problem size after the CPLEX presolve routine is listed under reduced model.}\label{table:comparison}
\begin{tabular}{lllll} \toprule
\textbf{Parameter}  & \multicolumn{2}{c}{\textbf{Initial model}}  & \multicolumn{2}{c}{\textbf{Reduced model}}\\ 
 & \text{Hourly} & \text{Cont.-time} & \text{Hourly} & \text{Cont.-time}\\  \midrule
\textbf{Binary variables} &       1,152 & 2,040 &      1,106  & 1,706      \\
\textbf{Continuous variables} &   2,474 & 8,954 &      2,179  & 7,371        \\
\textbf{Constraints}  &           2,706 &  16,962&       2,316 & 13,107       \\ \midrule
\textbf{Solution time [s]} &   &   & 2.2& 60.0  \\
\textbf{MIP gap [\%]} &  &    &        0.0 &0.28 \\ \bottomrule
\end{tabular}
\end{table}

\begin{figure}[ht]
    \centering 
    \includegraphics[width=\columnwidth]{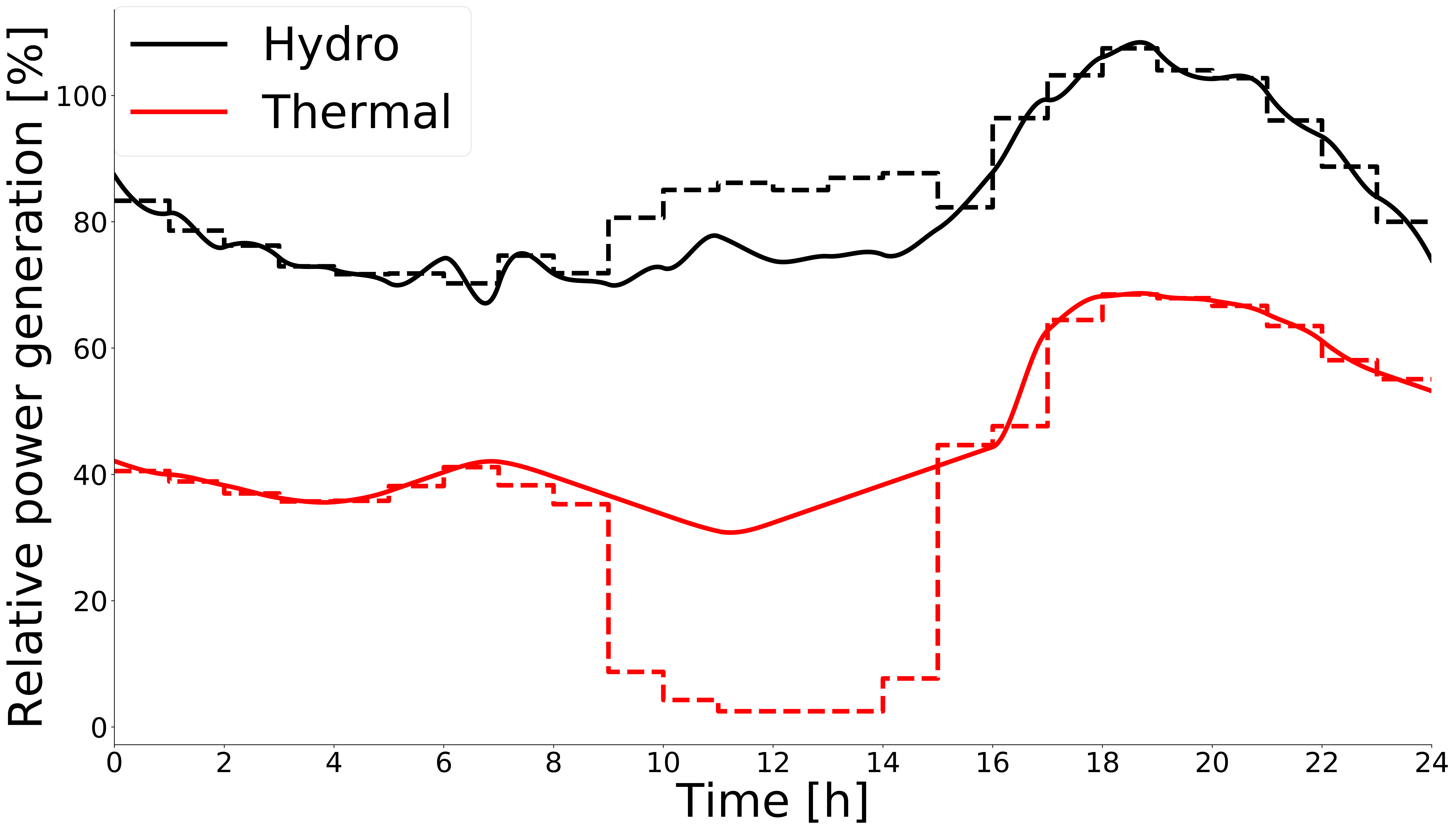}
\caption{The sum production in the thermal and hydropower areas relative to their respective load peaks in the hourly and continuous-time solution. The hourly and continuous-time solutions are shown as dashed and solid lines, respectively.}
\label{fig:production}
\end{figure}

The resulting sum production of hydropower and thermal generators are shown in \Cref{fig:production}. The figure shows that the hourly model overestimates the ramping capabilities of the thermal system during the extreme ramping events. Thermal production is shut down in the morning and turned back on in the afternoon, while the hydropower producers increase their production to cover the load in both areas in the meantime. This is not the case in the continuous-time model, as shutting down all thermal generators is either infeasible or very costly when following the net load during the ramping events. The cheapest and slowest thermal generator stays on for the whole 24 hours in the continuous-time model, contributing to the ramping in a modest way.

\begin{figure}[ht]
    \centering 
    \includegraphics[width=\columnwidth]{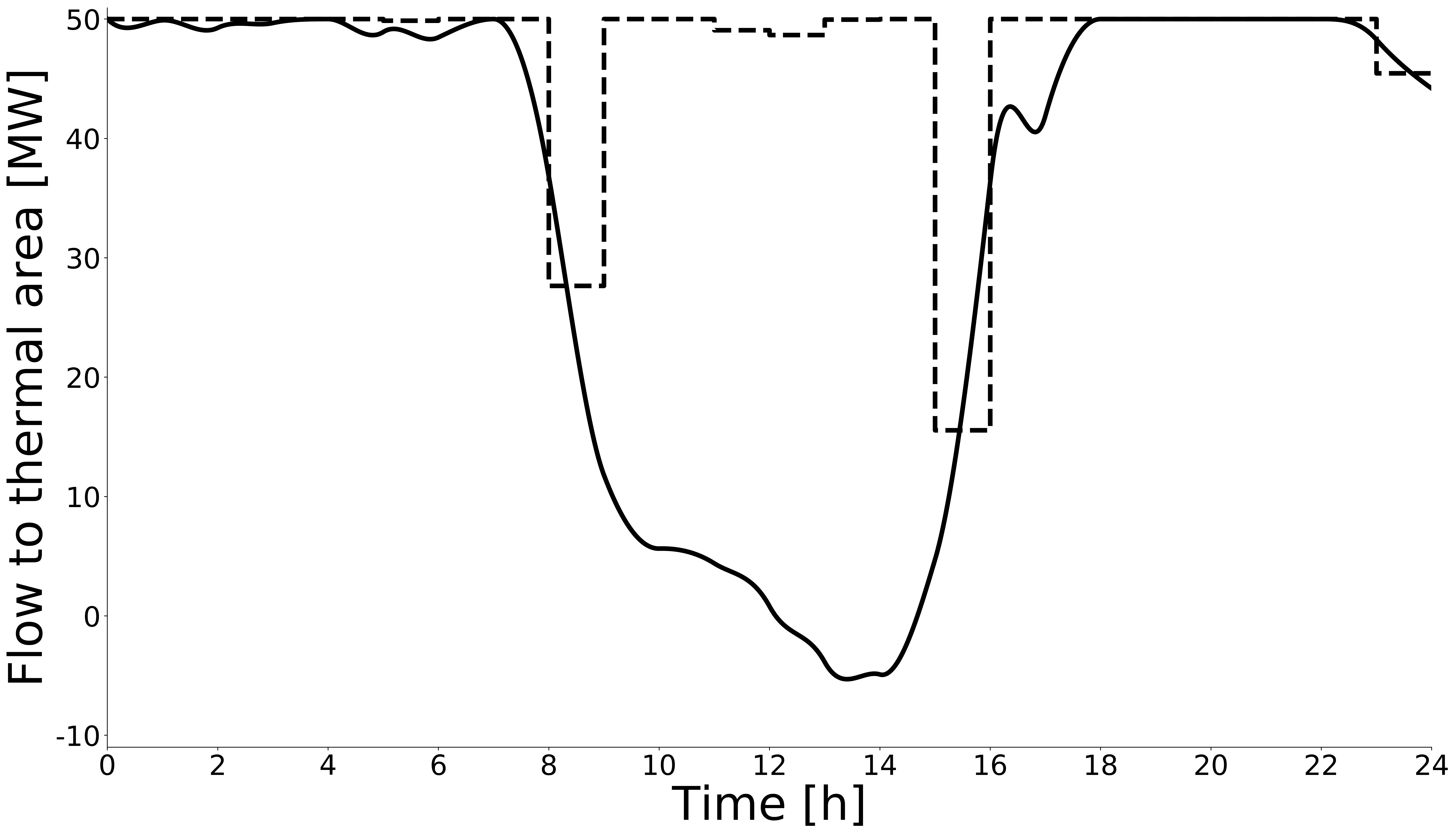}
\caption{The HVDC cable power flow from the hydropower area to the thermal area in the hourly (dashed line) and continuous-time (solid line) models. Negative values indicate flow in the opposite direction.}
\label{fig:flow}
\end{figure}

Most of the ramping is carried out by the hydropower system through the HVDC cable, which can be seen in \Cref{fig:flow}. The figure shows how the hydropower system is able to mitigate the ramping in net load in both directions while keeping the thermal generator online. The power flow is kept close to 50 MW throughout the day in the hourly model since the hydropower is generally cheaper than the thermal generators. However, two major changes in flow occur when the thermal generators are shut down and then started back up in the thermal system. This behaviour is undesirable, as it can increase the structural imbalances in the system \cite{statnett}.

\section{Conclusion}\label{conclusion_section}
Hydropower is considered an important balancing resource due to its flexibility. A continuous-time hydrothermal unit commitment model with HVDC cables was formulated in this paper to show how excessive ramping in the thermal system can be avoided by hydropower and active use of the HVDC cables. The structural imbalances in the system are reduced by 34\% in the continuous-time model compared to the hourly discrete-time model since sub-hourly effects are captured by the polynomial expansion. Several modelling issues related to incorporating hydropower into the continuous-time framework have been uncovered in the process. The linearization of the hydropower production curve requires binary variables to avoid unphysical uploading, and modelling the forbidden production zone requires the relaxation of the continuity constraints of the individual hydropower plants. The overall continuity of the model is still preserved on a system level, as the power balance forces the sum of hydropower production to be $C^1$ continuous. Investigating other potential modelling choices of the hydropower production curve, calculating system prices, and expanding the model to cover cross-zonal reserve capacity procurement are interesting avenues of further research. 



%

\bibliographystyle{IEEEtran}
\bibliography{IEEEabrv,refs}

\end{document}